\documentclass[12pt,intlim,righttag]{article}
\usepackage{amsmath,amsthm}
\usepackage{amssymb}
\usepackage{amsfonts}

\newcommand{\la}{\langle}
\newcommand{\ra}{\rangle}

\newcommand{\sg}{\sigma}

\newcommand{\lbl}{\label}

\newcommand\beq{\begin{equation}}
\newcommand\eeq{\end{equation}}

\newcommand{\beaa}{\begin{eqnarray*}}
\newcommand{\eeaa}{\end{eqnarray*}}

\theoremstyle{Theorem}

\theoremstyle{corollary}

\theoremstyle{remark}

\theoremstyle{definition}

\def\o{\omega}

\begin{document}
\title{Functional Equations for the Stochastic Exponential
 }

\author{B. Chikvinidze$^{1)}$, M. Mania$^{2)}$ and R. Tevzadze$^{3)}$}
\date{~}
\maketitle

\begin{center}
$^{1)}$ Georgian-American University and
Institute of Cybernetics of Georgian Technical Univercity,  Tbilisi, Georgia
\newline(e-mail: beso.chiqvinidze@gmail.com)
\\
$^{2)}$ Razmadze Mathematical Institute of Tbilisi State University  and
Georgian-American University, Tbilisi, Georgia,
\newline(e-mail: misha.mania@gmail.com)
\\
$^{3)}$ Georgian-American University and
Institute of Cybernetics of Georgian Technical Univercity,  Tbilisi, Georgia
\newline(e-mail: rtevzadze@gmail.com)
\end{center}

\begin{abstract}
{\bf Abstract.} We consider stochastic versions of the Cauchy exponential functional equation
 and give a martingale characterization of the general solution.
\end{abstract}

\noindent {\it 2010 Mathematics Subject Classification. 60G44, 60J65, 97I70}

\

\noindent {\it Keywords}: Stochastic Exponential, Martingales, Functional Equations, Non-anticipative functionals, Brownian Motion.

\section{Introduction}
For a continuous semimartingale $X=(X_t, t\ge0)$, with $X_0=0$,  the stochastic exponential is defined as
 \begin{equation}\label{st}
{\cal E}_t(X)= e^{X_t-\frac{1}{2}\la X\ra_t},\;\;\;t\ge0,
\end{equation}
where $\la X\ra$ is the square characteristic of the martingale part of $X$.

${\cal E}_t(X)$ is the unique solution of the linear stochastic differential equation
$$
Z_t=1+\int_0^t Z_s dX_s,\;\;\;t\ge0
$$
and in this sense $Z_t={\cal E}_t(X)$ is the stochastic analogue of the usual exponential function $f(x)=e^{cx}$, which is the unique
solution of the linear differential equation $f_x(x)=c f(x), f(0)=1.$

On the other hand, it is well known that the property
  \begin{equation}\label{exp}
e^{c(x+y)}=e^{cx}e^{cy},\;\;\;\text{for all}\;\;x,y\in R
\end{equation}
is a  characterizing property of exponential functions and in a very wide class of functions (e.g., in the class of measurable functions)  $f(x)=e^{cx}$, where $c\in R$ is some constant, is the general  solution
of the Cauchy exponential functional equation
\begin{equation}\label{fe}
f(x+y)=f(x) f(y)\;\;\;\text{for all}\;\;x,y\in R.
\end{equation}
Similar to (\ref{exp}) property of the stochastic exponential  (\ref{st}) is the equality
 \begin{equation}\label{st}
{\cal E}_t(X){\cal E}_t(Y)={\cal E}_t(X+Y+\la X, Y\ra),\;\;\;t\ge0,
\end{equation}
which is valid for any pair $X,Y$ of  semimartingales, where the additional term $\la X, Y\ra$ is the mutual characteristic of semimartingales $X$ and $Y$. For the properties of
stochastic exponential and for all unexplained notations from the martingale theory, we refer to \cite{Ja} or  \cite{LSH} .

By definition, the stochastic exponential ${\cal E}_t(X)$ of a continuous semimartingale $X$ can be expressed as a two-dimensional exponential function of a semimartingale and its square characteristic at time $t$

 \begin{equation}\label{st2}
{\cal E}_t(X)= f(\la X\ra_t, X_t),\;\;\;t\ge0,
\end{equation}
where
$f(u,v)=e^{v-\frac{1}{2}u},\;\;\;u\ge0, v\in R.$
 It follows from (\ref{st}), that this function satisfies the following functional equation
  \begin{equation}\label{fex}
f(\la X\ra_t, X_t)f(\la Y\ra_t, Y_t)=f(\la X+Y\ra_t, X_t+Y_t+\la X, Y\ra_t),
\end{equation}
valid for any  continuous semimartingales $X$ and $Y$ vanishing at $0$.
 It is easy to verify that (\ref{fex}) is also true for the function
\begin{equation}\label{fu}
f(u,v)=e^{cv-\frac{c}{2}u},\;\;\;u\ge0, v\in R,
\end{equation}
for any constant $c\in R.$

Denote by $\cal S$ the class of continuous semimartingales  and let ${\cal V}$ be  a subclass of $\cal S$, for which equation (\ref{fex}) is satisfied. Following \cite{K} we call the class ${\cal V}$ the
domain of validity of  equation (\ref{fex}).

One is led to the following question:
given a class of two-dimensional functions (e.g., measurable, continuous), how  small  a  class of
semimartingales $\cal V$ can we take, so that the solutions of (\ref{fex}) remain to be only  functions of the form (\ref{fu}).

We show in Theorem 1, that if equation (\ref{fex}) is satisfied for  the class of stochastic integrals $h\cdot W$  with respect to the given Brownian Motion $W$ and with
deterministic integrands $0\le h\le1$, then
the solution of equation (\ref{fex}) is of the form (\ref{fu}). In section 2 we prove this theorem by using the Cauchy exponential functional equation and the
corresponding  ''almost'' version.
In  section 3 we give a martingale characterization of equation (\ref{fex}), which gives also a probabilistic proof of  this theorem.

If the domain  of validity  of equation (\ref{fu}) consists only of a  single Brownian Motion $W$ (${\cal V}=\{W\}$),
which corresponds (if we take $X=Y=W$  in (\ref{fex})) to the equation
\begin{equation}\label{W1}
f^2(t, W_t)=f(4t, 2W_t+t),\;\;\;t\ge0,
\end{equation}
then there exists a
different from (\ref{fu}) continuous solution of equation (\ref{fex}), although any analytic solution of (\ref{W1}) is of the form (\ref{fu}). See Proposition 1 and a counterexample before
this proposition.

Note that the stochastic exponential ${\cal E}_t(X)$  transforms the class of continuous semimartingales into itself and this mapping is non-anticipative
in the following sense:
The mapping $F:R_+\times{\cal S}\rightarrow R$ is non-anticipative, if for any continuous semimartingales $X$ and $Y$ and $t\ge0$
$$
F(t, X)=F(t,Y), \;\;\;\text{when}\;\;\; X_s=Y_s\;\;\;\text{for all}\;\;\; s\le t.
$$

Therefore, it seems natural to consider   a functional equation for  stochastic exponents in terms of non-anticipative functionals
 \begin{equation}\label{fex2}
F(t, X)F(t, Y)=F(t, X+Y+\la X, Y\ra),
\end{equation}
 where by the property (\ref{st}) the stochastic exponential  ${\cal E}_t(X)$  satisfies this equation and  it represents  more general form of equation (\ref{fex}). But (\ref{fex2}) will be not a
 characterizing property of the
stochastic exponent, since there exists a whole class of solutions of (\ref{fex2})  which are not stochastic exponentials. E.g., if
\beq\lbl{ddol3}
F(t,X)=e^{\int_0^t{\cal K}(t,s)d(X_s-\frac12\la X\ra_s)},\;\;\;X\in\cal S,
\eeq
where $(k(t,s), s\ge0, t\ge0)$  is bounded, measurable deterministic function, then  $F(t,X)$  defined by (\ref{ddol3}) satisfies
 (\ref{fex2}), but such processes are not always stochastic exponentials, see Theorem 3 of section 4 and Remark 2 at the end of the paper. In Theorem 3
we prove that under some  restriction on the class of  non-anticipative functionals, the general solution of (\ref{fex2}) is of the form  (\ref{ddol3}).

\section{Functional equation for a function of a semimartingale and its square characteristic}

Let $W=(W_t, t\ge 0)$ be a standard Brownian Motion defined on a  complete probability space  $(\Omega, \cal F, P)$. Let
  $F=({{\cal F}}_t,t\ge0)$ be a filtration satisfying  the usual conditions of right-continuity and completeness. Let $F^W=({{\cal F}}^W_t,t\ge0)$ be the  filtration generated by the Brownian Motion $W$.

Let $\cal S$ (resp. $\cal M$) be the class of continuous semimartingales (martingales) vanishing at $0$.

Denote by ${\cal M}^W$ the class of continuous local martingales $M=(M_t, t\ge0)$ adapted to the filtration $F^W$  with $M_0=0$, i.e., this is the class of stochastic integrals  with respect to the Brownian Motion $W$.

Let ${\cal M}^W(I)$ be a sub-class of stochastic integrals $h\cdot W$ with respect to the Brownian Motion $W$ with
 integrands $h$, such that $h_u=I_{[s\le u\le t]}, 0\le s\le t$.

 Let $\left(f\left(u,v\right),u\ge0, v\in R\right)$ be a  function of two variables . We consider the  functional equation
  \begin{equation}\label{lob}
f(\la X\ra_t, X_t)f(\la Y\ra_t, Y_t)=f(\la X+Y\ra_t, X_t+Y_t+\la X, Y\ra_t),
\end{equation}
for any $X, Y\in {\cal V}$  and
$P$- a.e. for each $t\ge0$, where $\cal V$ is some class of  continuous semimartingales, the domain of validity of equation (\ref{lob}).

{\bf Theorem 1.} Let $\left(f\left(u,v\right),u\ge0, v\in R\right)$ be a  function of two variables. Then the following assertions are equivalent:

$a)$ The function $f$ is a continuous  strictly positive solution of the functional equation (\ref{lob}) with the domain of validity $\cal V=\cal S$.

$b)$ The function $f$ is a continuous  strictly positive solution of the functional equation (\ref{lob}) with the domain of validity ${\cal V}={\cal M}^W(I)$.

$c)$ The function $f$ is of the form
$$
f(u,v)=e^{cv - \frac{c}{2}u}\;\;\;\;\text{for some constant}\;\;\; c\in R.
$$

If we shall consider only measurable solutions, then the following two conditions will be equivalent:

$b')$ The function $f$ is a measurable  strictly positive solution of the functional equation (\ref{lob}) with the domain of validity ${\cal V}={\cal M}^W(I)$.

$c')$ The function $f=(f(u,v), u\ge0,v\in R)$ coincides  with the function
\begin{equation}\label{ec}
e^{cv - \frac{c}{2}u}\;\;\;\;\text{for some constant}\;\;\; c\in R
\end{equation}
almost everywhere with respect to the Lebesgue measure on $R_+\times R$.

\begin{proof} The  implication $a)\to b)$ is evident. Let us prove the implication $b)\to c)$.

It follows from equation (\ref{lob}) that for any bounded deterministic functions $h$ and $g$
 \begin{equation}\label{lob1}
f\Big(\int_0^th^2_udu, \int_0^th_udW_u\Big)
 f\Big(\int_0^tg^2_udu, \int_0^tg_udW_u\Big)=
\end{equation}
$$
=f\Big(\int_0^t(g_u+h_u)^2du, \int_0^t(g_u+h_u)dW_u+\int_0^t g_uh_udu\Big)
$$
$P$- a.e. for each $t\ge0$.

For any fixed pair $s\le t$ if we take $h_u=I_{(u<s)}$ and $g_u=I_{s\le u\le t)}$, from (\ref{lob1}) we obtain that $P$-a.s.
 \begin{equation}\label{lob2}
f(s, W_s)f(t-s, W_t-W_s)=f(t, W_t).
\end{equation}

From (\ref{lob2}) we have that
 \begin{equation}\label{lob21}
0=EI_{(f(s, W_s)f(t-s, W_t-W_s)\neq f(t, W_t))}=
\end{equation}
$$
=\int_R\int_R I_{(f(s, x)f(t-s, y)\neq f(t, x+y))}\rho(s,x)\rho(t-s,y-x)dxdy,
$$
where $\rho(s,x)=\frac{1}{\sqrt{2\pi s}}e^{\frac{x^2}{2s}}$.

Therefore, for any $s>0,t>0, s\le t$
 \begin{equation}\label{coshy2}
f(s, x)f(t-s, y)=f(t, x+y)\;\;\;\text{for all}\;\;\;x,y\in R
\end{equation}
almost surely with respect to the Lebesgue measure on $R^2$.

If $f(t,x)$ is continuous, then  (\ref{coshy2}) is fulfilled for all  $s>0,t>0, s\le t, x\in R, y\in R$ and $f$ satisfies the two-dimensional Cauchy exponential
functional equation. The general continuous solution of this equation is of the form (see, e.g., \cite{A})
\begin{equation}\label{W}
f(t, x)=\exp\{cx+bt\}\;\;\;\text{for some}\;\;\;b,c\in R.
\end{equation}
On the other hand inserting $X=W$ and $Y=W$  in (\ref{lob})  we have
\begin{equation}\label{form21}
f^2(t, W_t)=f(4t, 2W_t+t)\;\;\;t\ge0
\end{equation}
$P$-a.s. and substituting $f(t, x)=\exp\{cx+bt\}$ in (\ref{form21}) we obtain that
$$
\exp\{2cW_t+2bt\}=\exp\{2cW_t+ct+4bt\},
$$
which implies that $b=c/2$, hence $f$ is of the form $\exp\{cx-\frac{c}{2}t\}$.

$c)\to a)$ is  evident. Indeed, if $f(u,v)=\exp\{cv-\frac{c}{2}u\}$ for some $c\in R$, then for any $X,Y\in\cal S$
$$
f(\la X\ra_t, X_t)f(\la Y\ra_t, Y_t)=\exp\{cX_t-\frac{c}{2}\la X\ra_t\}\exp\{cY_t-\frac{c}{2}\la Y\ra_t\}=
$$
$$
=\exp\{c(X_t+Y_t)-\frac{c}{2}\la X+Y\ra_t +c\la X,Y\ra_t\}=f(\la X+Y\ra_t, X_t+Y_t+\la X, Y\ra_t)
$$
$P$-a.s. for any $t\ge0$.

The proof of the second part of Theorem 1 is similar if we use corresponding results on  $''$almost$''$ solutions of equation (\ref{coshy2}).

\
\\
$b')\to c')$.  It follows from results of \cite{DB} and \cite{J} that (\ref{coshy2}) implies
\begin{equation}\label{W2}
f(t, x)=\exp\{cx+bt\}\;\;\;\text{for some}\;\;\;b,c\in R,
\end{equation}
almost everywhere with respect to the Lebesgue measure on $R_+\times R$.

From (\ref{form21}) we have that
\begin{equation}\label{form22}
f^2(t, x)=f(4t, 2x+t)\;\;\;t\ge0
\end{equation}
almost everywhere with respect to the Lebesgue measure on $R_+\times R$. Therefore,
from (\ref{W2}) and (\ref{form22}), avoiding three null sets, we obtain that
$$
\exp\{2cx+2bt\}=\exp\{2cx+ct+4bt\},
$$
almost everywhere with respect to the Lebesgue measure on $R_+\times R$, which implies that $b=c/2$.

$c')\to b')$. Since for any $M\in{\cal M}^W(I)$ the random variable $M_t$ admits Gaussian distribution, it follows
from $c')$ that
\begin{equation}\label{f2}
f(\la M\ra_t, M_t)=\exp\{cM_t-\frac{c}{2}\la M\ra_t\}
\end{equation}
$P$-a.s. for any $t\ge0$, which will imply this implication, similarly to the proof of corresponding implication in first part of Theorem 1.

\end{proof}

{\bf Remark.} Note that, equality (\ref{f2}) can not be deduced from  condition $c')$ for a general continuous semimartingale $X$.

\

If the domain  of validity  of equation (\ref{lob}) consists only of a  single Brownian Motion $W$ (${\cal V}=\{W\}$), which corresponds to the equation
\begin{equation}\label{form33}
f^2(t, W_t)=f(4t, 2W_t+t),\;\;\;t\ge0, P-a.s.
\end{equation}
then there exists a
different from (\ref{fu}) continuous solution of equation (\ref{form33}).
Let $g(t,x)=\ln f(t,x)$. Then (\ref{form33}) is equivalent to equation
\begin{equation}\label{form44}
g(t, W_t)=\frac{1}{2}g(4t, 2W_t+t),\;\;\;t\ge0.
\end{equation}

Let $g(t,x)=|(2x-t)t|^\frac13$. It is evident, that this function satisfies equation (\ref{form44}), since
\[
g(4t,2x+t)=|(4x-2t)4t|^\frac13=2|(2x-t)t|^\frac13=2g(t,x).
\]
Therefore, $f(t,x)=\exp{|(2x-t)t|^\frac13}$ is a continuous solution of equation (\ref{form33}) different from (\ref{fu}).

$''$ It has been shown in various works that when some additional smoothness assumptions are imposed on set of functions of
classical functional equations,
then even if the domain of validity is quite small, the set of solution does not grow $''$ (see \cite{SH} for this comment and the references therein).
In the following proposition we show that if the domain of validity contains only a single Brownian Motion ($\cal V=W$), then
any analytic solution of (\ref{form33}) is of the form (\ref{fu}), i.e., the set of solutions is the same as for $\cal V=\cal S$.

{\bf Proposition 1.}
If $f$ is analytical  function satisfying (\ref{form33}), then $f$ is of the form  (\ref{fu}).

\begin{proof}
It follows from (\ref{form22}) and continuity of $f(t,x)$, that $f$ satisfies equation
\begin{equation}\label{form23}
f^2(t, x)=f(4t, 2x+t)\;\;\;\text{for all}\;\;\;t\ge0, x\in R.
\end{equation}
Taking logarithms of both sides of  (\ref{form23})  we obtain functional equation
\begin{equation}\label{form24}
h(t, x)=\frac{1}{2}h(4t, 2x+t)\;\;\;\text{for all}\;\;\;t\ge0, x\in R,
\end{equation}
for $h(t,x)=\ln f(t,x)$.

It is evident that $h(0,0)=0$. Differentiating equality (\ref{form24})  at $t$ we have
\begin{equation}\label{form34}
h_t(t,x)=2h_t(4t,2x+t)+\frac{1}{2}h_x(4t,2x+t).
 \end{equation}
This implies that  $h_t(0,0)=2h_t(0,0)+h_x(0,0)$ and hence
\begin{equation}\label{form25}
 h_t(0,0)=-\frac{1}{2}h_x(0,0).
 \end{equation}
 It easy to see that all higher order derivatives of $h$ at $t=0,x=0$ are equal to zero.

 Indeed, it follows from (\ref{form24}) that for any $m\ge1$ and $n\ge1$
 $$
 \frac{\partial^{m+n}}{\partial t^m\partial x^n}h(t,x)= \frac{\partial^{m+n-1}}{\partial t^m\partial x^{n-1}}h_x(t,x)=
  \frac{\partial^{m+n-1}}{\partial t^m\partial x^{n-1}}h_x(4t,2x+t)=
 $$
 $$
=\frac{1}{2}  \frac{\partial^{m+n}}{\partial t^m\partial x^{n}}\big(h(4t,2x+t)\big)_x=\frac{1}{2}  \frac{\partial^{m+n}}{\partial t^m\partial x^{n}}h(4t,2x+t),
 $$
which implies that
$$
 \frac{\partial^{m+n}}{\partial t^m\partial x^n}h(0,0)=\frac{1}{2} \frac{\partial^{m+n}}{\partial t^m\partial x^n}h(0,0)\;\;\;\text{and}\;\;\;
  \frac{\partial^{m+n}}{\partial t^m\partial x^n}h(0,0)=0.
  $$
If $n=0$ then it follows from (\ref{form24}) and (\ref{form34}) that
$$
 \frac{\partial^{m}}{\partial t^m}h(t,x)= \frac{\partial^{m-1}}{\partial t^{m-1}}h_t(t,x)=
 \frac{\partial^{m-1}}{\partial t^{m-1}}\big(2h_t(4t,2x+t)+\frac{1}{2}h_x(4t,2x+t)\big)=
 $$
$$
=\frac{1}{2} \frac{\partial^{m}}{\partial t^{m}}h(4t,2x+t)+\frac{1}{4} \frac{\partial^{m}}{\partial t^{m-1}\partial x}h(4t,2x+t).
$$
Since   $\frac{\partial^{m+n}}{\partial t^m\partial x^n}h(0,0)=0$ for any $n\ge1$, this implies that
$$
\frac{\partial^m}{\partial t^m}h(0,0)=\frac{1}{2}\frac{\partial^m}{\partial t^m}h(0,0)
$$
and this derivative at point $(0,0)$ is also equal to zero.

 Therefore, since $f$ is analytic, we obtain  from (\ref{form25}) that
$$
 h(t,x)=h_t(0,0)t+f_x(0,0)x=h_x(0,0)x-\frac{h_x(0,0)}{2}t,
 $$
which means that $h$ is of the form $h(t,x)=cx-\frac{c}{2}t$ for some $c\in R$.
\end{proof}

\section{A martingale approach}

In this section we give a martingale characterization of equation (\ref{lob}), which gives also a probabilistic proof of Theorem 1.

Let first show that $f(t,W_t)$  is integrable at any power.

{\bf Lemma 1.}  If  $\left(f\left(u,v\right),u\ge0, v\in R\right)$ is a  measurable  strictly positive solution of the functional equation (\ref{lob}), then
$$
E|f(t,W_t)|^p<\infty
$$
for any $t\ge0$ and $p\in R$.

\begin{proof}
To show that $f\left(t, W_{t}\right)$ is integrable for any $t\geq 0$  we shall use the idea from \cite{SM} on application of the Bernstein theorem.

Substituting $s=t/2$ in (\ref{lob2}) and taking logarithms we have
 \begin{equation}\label{t2}
\ln f(t/2, W_{t/2})+\ln f(t/2, W_t-W_{t/2})= \ln f(t, W_t) \;\;\;\;P-a.s.
\end{equation}
Let
$$
X=\ln f(t/2, W_{t/2})\;\;\;\text{and}\;\;\;Y=\ln f(t/2, W_t-W_{t/2}).
$$
Then from (\ref{t2})
 \begin{equation}\label{t3}
X+Y= \ln f(t, W_t)\;\;\;\;P-a.s.
\end{equation}
and it is easy to see that
$$
X-Y=\ln f(t/2, W_{t/2})-\ln f(t/2, W_t-W_{t/2})= \ln f(t, 2W_{t/2}- W_t)\;\;\;\;P-a.s.
$$

Since $E(2W_{t/2}- W_t)W_t=0$, the random variables $2W_{t/2}- W_t$ and $W_t$ are  independent and, hence, the random
variables $X+Y$ and $X-Y$ will also be independent.
Therefore, Bernstein's theorem (\cite{BE}, see also \cite{Q}) implies that $X$ and $Y$ are distributed normally. Therefore,  $\ln f(t, W_t)$ is also normally distributed. Hence
$f(t,W_t)$ is integrable at any power
$$
E|f(t, W_{t})|^p<\infty,
$$
as a random variable having log-normal distribution.
\end{proof}

{\bf Proposition 2.} If  $\left(f\left(u,v\right),u\ge0, v\in R\right)$ is a  measurable  strictly positive solution of the functional equation (\ref{lob}), then

$i)$ the  processes
$$
N_t(1)=\frac{f(t, W_t)}{Ef(t, W_t)},\;\;\;\text{and}\;\;\; N_t(2)=\frac{f^2(t, W_t)}{Ef^2(t, W_t)},\;\;t\ge0,
$$
are   strictly positive martingales.

$ii)$ There exist constants $\lambda_1\in R$ and $\lambda_2\in R$ such that the processes
$$
\frac{f(t, W_t)}{e^{\lambda_1 t}},\;\;\;\frac{f^2(t, W_t)}{e^{\lambda_2 t}}\;\;\;\text{and}\;\;\;\frac{f(4t, 2W_t+t)}{e^{\lambda_2 t}}    \;\;\; t\ge0,
$$
are strictly positive martingales.  Besides
\begin{equation}\label{12}
\lambda_1\le \frac{1}{2}\lambda_2.
\end{equation}

\begin{proof}

$i)$ Taking expectations in equality (\ref{lob2}), since the random variables $W_t-W_s$  and $W_s$ are  independent, we have that
\begin{equation}\label{lob4}
Ef(t-s, W_t-W_s)=\frac{E f(t,W_t)}{Ef(s, W_s)}.
\end{equation}

If we take the conditional expectations in the same equality (\ref{lob2}),  having in mind that $W_t-W_s$ is independent of ${\cal F}_s^W$, we obtain
\begin{equation}\label{lob5}
E(f(t, W_t)| {{\cal F}_s^W}) =f(s, W_s)Ef(t-s, W_t-W_s)\;\;\;\;P-a.s..
\end{equation}
Therefore,  substituting the expression of $Ef(t-s, W_t-W_s)$ from  (\ref{lob4}) into  (\ref{lob5}) we get the martingale equality
$$
E\Big(\frac{f(t, W_t)}{Ef(t,  W_t)}\Big|{\cal F}_s\Big)=\frac{f(s,W_s)}{Ef(s,W_s)},  \;\;\;\;P-\text{a.s.}.
$$
The martingale property of the process $N_t(2)$ is proved similarly, if we use equality

 \begin{equation}\label{lobp}
f^2(s, W_s)f^2(t-s, W_t-W_s)=f^2(t, W_t)\;\;\;\;P-a.s.
\end{equation}
instead of equality (\ref{lob2}).

$ii)$ Let $g(t)\equiv E f(t,W_t)$. Then $Ef(t-s, W_t-W_s)=g(t-s)$ and it follows from (\ref{lob4}) that $(g(t), t\ge0)$ satisfies the  Cauchy
exponential functional equation
$$
g(s)g(t-s)=g(t),\;\;\;s\ge0, t\ge0, s\le t
$$
a general measurable solution of which is of the form $e^{\lambda t}$ for some $\lambda\in R$ (see, e.g., \cite{A}). Therefore the proof of $ii)$ follows from assertion $i)$.

Since $f^2(t, W_t)e^{-\lambda_2 t}$ is a martingale, the equality (\ref{form21}) implies that the process $f(4t, 2W_t+t)e^{-\lambda_2 t}$ is also  a martingale.

The inequality (\ref{12}) follows from the martingale property of $N(1)$ and $N(2)$ and from the H\^older inequality, since
$$
1=E\frac{f(t, W_t)}{e^{\lambda_1 t}}\le E^{1/2}\frac{f^2(t, W_t)}{e^{2\lambda_1 t}}=
$$
$$
=E^{1/2}\frac{f^2(t, W_t)}{e^{\lambda_2 t}}\cdot\frac{e^{\lambda_2 t}}{e^{2\lambda_1 t}}= e^{\big(\frac{\lambda_2}{2}-\lambda_1\big) t}.
$$
\end{proof}

{\bf Theorem 2.} Each assertion $a), b), c)$ of Theorem 1 is equivalent to

$d)$  $\left(f\left(u,v\right),u\ge0, v\in R\right)$ is a  continuous  strictly positive function such that  $f(0,0)=1$ and the processes
$$
\frac{f(t, W_t)}{e^{\lambda_1 t}},\;\;\;\frac{f^2(t, W_t)}{e^{\lambda_2 t}},\;\;\;\text{and}\;\;\;\frac{f(4t, 2W_t+t)}{e^{\lambda_2 t}},\;\;\; t\ge0,
$$
are  martingales for some constants $\lambda_1\in R$ and $\lambda_2\in R$.

\begin{proof} The implication $b)\to d)$ follows from Proposition 2. Let us show that if $d)$ is satisfied then the function $f(t,x)$ is of the form (\ref{ec}).

 Let us define two functions:
$$
U(t,x)=E\big[ f(T,W_T)e^{-\lambda_1 T} \big| W_t = x \big]\;\;\text{and}\;\;V(t,x)=E\big[ f^2(T,W_T) e^{-\lambda_2 T} \big| W_t = x \big].
$$
Since $U(t,x)$ and $V(t,x)$ are positive, they will be of the class $C^{1.2}$ on $(0,T)\times R$ and  satisfy the ''backward'' heat equation (see, e.g. \cite{KAR} page 257)
\begin{equation}\label{kol}
\frac{\partial Y}{\partial t}+\frac{1}{2}\frac{\partial^2Y}{\partial x^2}=0,\;\;\;\;0<t<T, x\in R.
\end{equation}
Since the processes $f(t,W_t)e^{-\lambda_1 t}$ and $f^2(t,W_t)e^{-\lambda_2t}$ are martingales, by the Markov property of $W$ we shall have that
$$
U(t,W_t)={f(t, W_t)}{e^{-\lambda_1 t}}\;\;\;\text{and}\;\;\; V(t,W_t)=f^2(t,W_t)e^{-\lambda_2 t}
$$
and  hence,   by continuity of the function $f$
$$
U(t,x)={f(t, x)}{e^{-\lambda_1 t}}\;\;\;\text{and}\;\;\; V(t,x)={f^2(t, x)}{e^{-\lambda_2 t}},\;\;\;\text{for all}\;\;\; t\ge0, x\in R.
$$

Since
$$
U_t(t,x)=f_t(t,x)e^{-\lambda_1 t}-\lambda_1 f(t,x)e^{-\lambda_1 t}; \;\;\; U_{xx}(t,x)=\frac{1}{2}f_{xx}(t,x)e^{-\lambda_1 t}
$$
and $U(t,x)$ satisfies equation (\ref{kol}), we shall have that  $f$ also belongs to the class $C^{1.2}$ on $(0,T)\times R$ and  satisfies equation
\begin{equation}\label{eq1}
f_t(t,x)+\frac{1}{2}f_{xx}(t,x)-\lambda_1 f(t,x)=0.
\end{equation}

Similarly, since
$$
V_t(t,x)=2f(t,x)f_t(t,x)e^{-\lambda_2t} - \lambda_2 f^2(t,x)e^{-\lambda_2t};
$$
$$
V_{xx}(t,x)= 2(f_x(t,x))^2 e^{-\lambda_2 t} + 2f(t,x)f_{xx}(t,x)e^{-\lambda_2 t}
$$
and $V(t,x)$ satisfies equation (\ref{kol}), we obtain that the function $f$ satisfies also the following PDE:
$$
2f(t,x)f_t(t,x)e^{-\lambda_2t} - \lambda_2 f^2(t,x)e^{-\lambda_2t} +
$$
$$
\frac{1}{2}\Big[ 2(f_x(t,x))^2 e^{-\lambda_2 t} + 2f(t,x)f_{xx}(t,x)e^{-\lambda_2 t}\Big] = 0,
$$
which after simplifying takes the form:
\begin{equation}\label{eq2}
2f(t,x)f_t(t,x) - \lambda_2 f^2(t,x) + (f_x(t,x))^2 + f(t,x)f_{xx}(t,x) = 0.
\end{equation}

From (\ref{eq1}) we have that $f_{xx}(t,x)=2\lambda_1 f(t,x)-2f_t(t,x)$ and  substituting the right side of this equality instead of $f_{xx}$ in (\ref{eq2})
we get that
\begin{equation}\label{eq3}
f^2_x(t,x)=(\lambda_2 - 2\lambda_1)f^2(t,x),
\end{equation}
\\
where $\lambda_2- 2 \lambda_1\ge0$ according to Proposition 2.

To find the general solution of (\ref{eq3}) we consider two cases:
\\
I: $\; f_x(t,x)=\sqrt{\lambda_2 - 2\lambda_1}f(t,x)$. It is evident that $f(t,x)=\phi (t)e^{\sqrt{\lambda_2 - 2\lambda_1}x}$, where $\phi$ is
some function of $t$. It follows from Proposition 2 that $Ef(t,W_t)=e^{\lambda_1 t}$. So
$e^{\lambda_1 t}=\phi(t)Ee^{\sqrt{\lambda_2 - 2\lambda_1}W_t}=\phi(t)e^{\frac{\lambda_2-2\lambda_1}{2}t}$, which implies that
$\phi(t)=e^{\frac{4\lambda_1-\lambda_2}{2}t}$ and finally
$f(t,x)=e^{cx+bt}$, where $c=\sqrt{\lambda_2 - 2\lambda_1}$ and
$b=\frac{4\lambda_1-\lambda_2}{2}$.
\\
II: $\; f_x(t,x)=-\sqrt{\lambda_2 - 2\lambda_1}f(t,x)$.  By the same manner as in the case I, we obtain that
$f(t,x)=e^{cx+bt}$, where $c=-\sqrt{\lambda_2 - 2\lambda_1}$ and
$b=\frac{4\lambda_1-\lambda_2}{2}$.

After that we can use equality $Ef^2(t,W_t)=Ef(4t,2W_t+t)$ to get equality $b=-c/2$ and the representation $f(t,x)=e^{cx-\frac{c}{2}t}$.
\end{proof}

\section{A functional equation for non-anticipative functionals}

The mapping $h:[0,T]\times{C[0,T]}\rightarrow R$ is non-anticipative, if for any $\omega, \omega'\in C[0,T]$ and $t\in [0,T]$
$$
h(t, \omega)=h(t,\omega'), \;\;\;\text{when}\;\;\; \omega_s=\omega'_s\;\;\;\text{for all}\;\;\; s\le t.
$$

Consider the class of  functions $F:[0,T]\times{\cal{S}}\rightarrow R$ defined by
\beaa
{\cal C}=\{F;\;F(t,X)=e^{h(t,X-\frac12\la X\ra)},X\in{\cal S},\\
\text{for some continuous, non-anticipative functional}\; h(t,\o)\}.\eeaa

{\bf Theorem 3.} The general solution of the functional  equation
 \begin{equation}\label{fex3}
F(t, X)F(t, Y)=F(t, X+Y+\la X, Y\ra),\;\;\;\text{for all}\;\;\;X, Y\in\cal S,
\end{equation}
in the class ${\cal C}$, is of the form
\beq\lbl{ddol2}
F(t,X)=e^{\int_0^t{\cal K}(t,s)d(X_s-\frac12\la X\ra_s)},\;\;\;X\in\cal S,
\eeq
where $({\cal K}(t,s), s\ge0, t\ge0)$  is  a deterministic function with  ${\cal K}(t,s)=0,t\le s$, such that
\begin{description}
  \item[i)] ${\cal K}(t,\cdot)\;\;\text{ is  cadlag and has a finite variation,  for each}\;\;t\in[0,T],$
  \item[ii)] $ \int_0^T\o_s{\cal K}(\cdot,ds)\;\;\;\text{ is continuous, for each}\;\;\o\in C[0,T].$
\end{description}

\begin{proof}
Let $F$ be a solution of equation (\ref{fex3}) from the class $\cal C$. Then it follows from (\ref{fex3}) and  from the definition of the class $\cal C$ that
\beaa
h(t,X)+h(t,Y)=\ln F(t,X+\frac12\la X\ra)+\ln F(t,Y+\frac12\la Y\ra)\\
=\ln F(t,X+Y+\frac12(\la X\ra+\la Y\ra)+\la X,Y\ra)\\
=\ln F(t,X+Y+\frac12(\la X+Y\ra))=h(t,X+Y)
\eeaa
for any $X, Y\in {\cal S}$.

Since a deterministic function is a semimartingale if and only if it is of finite variation and the functions of finite variations
are dense in $C$,  it follows from the continuity of $h$ that
 \begin{equation}\label{c0t}
 h(t,\omega+\omega')={h}(t,\omega)+{h}(t,\omega') \;\;\text{for all}\;\;\omega, \omega'\in{C}.
\end{equation}
By the Riesz theorem (see, e.g. \cite{DSc}) for
each $t$ there exists a cadlag function  $G(t,\cdot)$ of finite variation, such that $G(t,s)=G(t,T),s\ge t$ and
\begin{equation}\label{var}
h(t,\o)=\int_0^T\o_sG(t,ds).
\end{equation}
Using integration by part formula we get from (\ref{var}) that
\beaa
h(t,X)=\int_0^tX_sG(t,ds)=X_tG(t,t)-\int_0^tG(t,s)dX_s\\
=\int_0^t(G(t,t)-G(t,s))dX_s=\int_0^t{\cal K}(t,s)dX_s,
\eeaa
where ${\cal K}(s,t)=G(t,t)-G(t,s)$ is of finite variation for each $t$ and by continuity of $h$ the condition $ii)$ is also satisfied.

Now let us show that the non-anticipative functional
$F(t,X)=e^{\int_0^t{\cal K}(t,s)d(X_s-\frac12\la X\ra_s)}$ satisfies equation (\ref{fex3}). Let first show that $F$ belongs to the class $\cal C$.

Consider the  non-anticipative functional  $h(t,\omega)=\int_0^T\o_s{\cal K}(\cdot,ds)$, where $\cal K$ satisfies conditions $i)-ii)$ of the theorem. By $ii)$ this functional is continuous  and the integration by part
formula gives the equality
 \begin{equation}\label{part}
h(t,X)=-\int_0^tX_s{\cal K}(t,ds)=\int_0^t{\cal K}(t,s)dX_s,\;X\in{\cal S}.
\end{equation}
Therefore, (\ref{part}) implies that
\beq\lbl{ddol}
F(t,X)=e^{\int_0^t{\cal K}(t,s)d(X_s-\frac12\la X\ra_s)}=e^{h(t,X- \frac12\la X\ra)},\;\;\;X\in\cal S,
\eeq
which means that $F\in\cal C$.

It is evident that $F(t,X)$  defined by (\ref{ddol2}) satisfies (\ref{fex3}), since for any continuous semimartingales $X$ and $Y$
$$
F(t, X)F(t, Y)=e^{\int_0^tk(t,s)d(X_s+Y_s-\frac12\la X\ra_s-\frac12\la Y\ra_s)}=
$$
$$
=e^{\int_0^tk(t,s)d(X_s+Y_s+\la X,Y\ra_s-\frac{1}{2}\la X+Y\ra_s)}=F(t, X+Y+\la X, Y\ra).
$$
\end{proof}

{\bf Remark 1.} If ${\cal K}(t,s)$, ${\cal K}(t,s)=0,t\le s$ is a measurable function, such that
$\int_0^T|{\cal K}(t,s)-{\cal K}(t',s)|^2ds\le const |t-t'|^2$,
then by the Kolmogorov theorem $\int_0^t{\cal K}(t,s)dX_s$ has a continuous modification, for each $X_t=X_0+\int_0^t\sg_sdW_s+\int_0^tb_sds$ with bounded, predictable $b,\sg$.

{\bf Remark 2.} Note that $F(t,X)$ is not always a stochastic exponential. For example, if we take ${\cal K}(t,s)= v(t)$,   where $v$ is  not of finite variation, then $F(t,X)=e^{v(t)(X_t-\frac12\la X\ra_t)}$
will be not a stochastic exponential for any $X\in{\cal S}$. E.g., for $X=W$ the process
$v(t)(W_t-\frac{1}{2}t)$ will be not a semimartingale and, hence $F(t, W)$  will not be a stochastic exponential.

{}

\end{document}